\documentclass[12pt,letter]{amsart}
\usepackage{amscd, ytableau}
\usepackage{amssymb,tikz-cd}
\usepackage[centering,text={15.5cm,23cm}]{geometry} 

\usepackage{amsfonts,amssymb, amscd, latexsym, graphicx, psfrag, color,float}

\usepackage{caption}
\usepackage{comment}
\usepackage{wrapfig}
\usepackage[all]{xy}
\usepackage{mathrsfs}
\usepackage{mathtools}
\usepackage{marvosym}
\usepackage{stmaryrd}
\usepackage{srcltx}
\usepackage{hyperref}
\usepackage{upgreek, MnSymbol} 
\usepackage[inline]{enumitem} 

\usepackage{braket}

\DeclareFontFamily{U}{mathx}{}
\DeclareFontShape{U}{mathx}{m}{n}{ <-> mathx10 }{}
\DeclareSymbolFont{mathx}{U}{mathx}{m}{n}
\DeclareFontSubstitution{U}{mathx}{m}{n}
\DeclareMathAccent{\widecheck}{0}{mathx}{"71}

\definecolor{mygray}{gray}{0.75} 

\definecolor{shadecolor}{rgb}{1,0.9,0.7}

\newtheorem{theorem}{Theorem}[section]
\newtheorem{lemma}[theorem]{Lemma}
\newtheorem{lemma-definition}[theorem]{Lemma-Definition}
\newtheorem{proposition}[theorem]{Proposition}
\newtheorem{corollary}[theorem]{Corollary}

\newtheorem{assumption}[theorem]{Assumption}
\newtheorem*{theoremx}{Theorem}
\newtheorem*{corollaryx}{Corollary}

\theoremstyle{definition}

\newtheorem{definition}[theorem]{Definition}

\newtheorem{example}[theorem]{Example}

\theoremstyle{remark}
\newtheorem{remark}[theorem]{Remark}

\numberwithin{equation}{section}
\numberwithin{figure}{section}

\newcommand{\bP}{\mathbb{P}}

\newcommand{\cO}{\mathcal{O}}

\newcommand{\barM}{\overline{\mathcal{M}}}
\newcommand{\barN}{\overline{N}}

\newcommand {\red}[1]{{\color{red} #1}}

\def\mydate{\ifcase\month \or January\or February\or March\or
April\or May\or June\or July\or August\or September\or October\or 
November\or December\fi \space\number\day,\space\number\year}

\makeatletter
\let\oldcite\cite
\renewcommand{\cite}{\@ifnextchar[{\@newcite}{\oldcite}}
\def\@newcite[#1]#2{\oldcite{#2},\,#1}
\makeatother

\begin{document}

\title[Tropicalization of super Gromov-Witten invariants]
{Tropicalization of super Gromov-Witten invariants}

\author{Artan Sheshmani \textsuperscript{1,2}}
\email{artan@mit.edu}

\author{Shing-Tung Yau \textsuperscript{1,3}} 
\email{styau@tsinghua.edu.cn}

\author{Benjamin Zhou \textsuperscript{3}}
\email{byzhou@mail.tsinghua.edu.cn}

\address{\textsuperscript{1}Beijing Institute of Mathematical Sciences and Applications, No. 544, Hefangkou Village, Huaibei Town, Huairou District, Beijing 101408}
\address{\textsuperscript{2}Massachusetts Institute of Technology, IAiFi Institute, 77 Massachusetts Ave, 26-555. Cambridge, MA 02139} 
\address{\textsuperscript{3}Yau Mathematical Sciences Center, Tsinghua University, Haidian District,
Beijing, China}

\begin{abstract}
We show that genus-0, $n$ Neveu-Schwarz marked, super Gromov-Witten invariants of a convex, toric variety $X$ can be defined and computed using tropical geometry. When $X$ is a point, the tropical, super Gromov-Witten invariants of $X$ are descendant invariants on the moduli space of tropical curves. When $X$ is a general convex, toric variety, we define a procedure that computes the tropical, inverse Euler class of the SUSY normal bundle $\barN_{n, \beta} \rightarrow \barM_{0,n}(X, \beta)$, under the assumption that $\barN_{n, \beta}$ is in some sense locally tropicalizable. We define the tropical, genus-0, $n$ Neveu-Schwarz marked, super Gromov-Witten invariants of $X$, and show that the definition recovers the tropical, super Gromov-Witten invariants of a point. We compute a tropical, genus-0, super Gromov-Witten invariant of $\bP^1$.
\end{abstract}

\maketitle
\setcounter{tocdepth}{1}
\tableofcontents

 
\section{Introduction}

Tropical geometry has been quite effective in translating questions of enumerative geometry to combinatorial ones. There have been theorems showing equivalences between counts of tropical curves and Gromov-Witten invariants of a target space $X$. When $X$ is a toric surface, Mikhalkin showed that genus-0, Gromov-Witten invariants of $X$ are counts of tropical curves in $\mathbb{R}^2$ \cite{Mik}. His work was generalized to higher dimensional toric varieties in \cite{NS}. Bousseau showed that $q$-refined counts of tropical curves in $\mathbb{R}^2$ compute $g \geq 0$, logarithmic Gromov-Witten invariants of $X$ \cite{Bou}. Tropical techniques are used to prove a degeneration formula for logarithmic Gromov-Witten invariants \cite{ACGS} \cite{KLR}, and to establish mirror symmetry for log Calabi-Yau surfaces \cite{GHK} \cite{GS16}. There is also work involving analytic geometry in tropicalization of the moduli of curves \cite{ACP} and of the moduli of stable maps \cite{Yu}. Tropical geometry has displayed a broad and effective role in enumerative geometry and mirror symmetry.

Suppose $X$ is a convex, toric variety with $\dim X = k$. Let $\barM_{0,n}(X, \beta)$ be the moduli space of genus-0, $n$-marked stable maps to $X$ in curve class $\beta \in H_2(X, \mathbb{Z})$. Gromov-Witten invariants are enumerative invariants of $X$ defined by integrating over the Deligne-Mumford stack $\barM_{0,n}(X, \beta)$ \cite{FP}. The invariants are virtual counts of stable maps $f:\Sigma \rightarrow X$ where $\Sigma$ is a possibly nodal, genus-$g$ Riemann surface. If $X$ is a symplectic manifold with a tame, almost complex structure $J$, Gromov-Witten invariants can also be defined using symplectic methods as counts of $J-$holomorphic curves to $X$. The definitions of Gromov-Witten invariants for complex projective varieties using algebraic geometry or symplectic geometry are equivalent \cite{Sie1} \cite{Sie2}.

It is natural to consider an extension of Gromov-Witten theory to maps from super Riemann surfaces. Supergeometry and super Riemann surfaces were studied in  \cite{LR} \cite{Wit}. The moduli space of stable super Riemann surfaces is a smooth Deligne-Mumford stack \cite{FKP}, and some of its intersection numbers have been computed \cite{Nor}. Recently, Keßler, Sheshmani, and Yau developed the genus-0, super Gromov-Witten theory of $X$. In \cite{KSY1}, they construct $\barM^{super}_{0,n}(X, \beta)$, the moduli space of genus-0, $n$ Neveu-Schwarz marked, super $J$-holomorphic maps to a convex, toric variety $X$ in curve class $\beta \in H_2(X, \mathbb{Z})$. In \cite{KSY2}, they show that $\barM_{0,n}(X, \beta)$ embeds as a $\mathbb{C}^*$-fixed locus into $\barM^{super}_{0,n}(X, \beta)$,

\begin{center}
\begin{tikzcd}
\label{diag:intro_SUSY}
\barN_{n, \beta}  \arrow[d] \\
\barM_{0,n}(X, \beta)  \arrow[hook, r] & \barM^{super}_{0,n}(X, \beta)
\end{tikzcd}
\end{center}
where $\barN_{n, \beta}$ is the \textit{SUSY normal bundle}, defined in Section \ref{sec:SUSY_bundles}. By $\mathbb{C}^*$-localization, the genus-0, $n$ Neveu-Schwarz marked \footnote{\label{fn:NS}Neveu-Schwarz marked points can be equivalently thought of as smooth marked points by \cite{FKP}, Section 2.5.}, super Gromov-Witten invariants of $X$ are obtained by integration over $\barM_{0,n}(X, \beta)$ of incidence conditions and the inverse equivariant Euler class of $\barN_{n, \beta}$ with respect to $\mathbb{C}^*$. The authors subsequently apply localization to $\barM_{0,n}(X, \beta)$ with the $(\mathbb{C}^*)^k$-torus action of $X$ to compute super Gromov-Witten invariants of $\bP^n$ \cite{KSY3}.

\subsection{Main results}
In this paper, we use tropical geometric methods to define and compute tropical, genus-0, super Gromov-Witten invariants of a convex, toric variety $X$. In Section \ref{sec:SUSY_SGW}, we recall the SUSY normal bundle, and review the definition of super Gromov-Witten invariants \cite{KSY3}. We prove facts about the SUSY normal bundle $\barN_{n, \beta} \rightarrow \barM_{0,n}(X, \beta)$, including Proposition \ref{prop:vector_bundle}, which shows that the SUSY normal bundle $\barN_{n, \beta}$ is a vector bundle of rank $c_1(\beta) + n - 2$. In Section \ref{sec:SGW_point}, we review descendant invariants on the moduli of curves $\barM_{0,n}$, as well as tropical descendant invariants on the moduli space of tropical curves $\barM^{trop}_{0,n}$. In Definition \ref{def:super_dual_graphs}, we define an enlarged moduli space $\overline{T}_{0,n}$, which takes into account powers of $\psi_i$-classes on a stable curve. In Definition \ref{def:SGW_point}, we define the tropical, genus-0, $n$-marked super Gromov-Witten invariant of a point, $TSGW(pt)$. Using results of \cite{KM}, we show that, 

\begin{theoremx}[Theorem \ref{thm:TSGW_point}]
    Let $X = pt$. We have,
 
    \[
    SGW_{0,n}(pt) = TSGW_{0,n}(pt)
    \]
    where $SGW_{0,n}(pt)$ is given in Proposition \ref{prop:SGW_point_descendant}, and $TSGW_{0,n}(pt)$ is defined in Definition \ref{def:tropical_SGW}.
\end{theoremx}
In Section \ref{sec:SGW_toric}, we define the tropical, genus-0, super Gromov-Witten invariants of a general convex, toric variety $X$. Let $\Sigma_X$ be the fan of $X$. Let $TSM_{0,n}(\Sigma_X, \beta)$ be the extended cone complex of tropical stable maps obtained from $\barM_{0,n}(X, \beta)$ \cite{Ran}. Let $K:= \mathbb{C}^*$ act on $\barM_{0,n}(X, \beta)$ as the identity, and let $\kappa: K \rightarrow \mathbb{C}^*$ be the identity character. In Theorem \ref{thm:tropical_Euler_class}, we define a procedure that computes the tropical, inverse Euler class of the SUSY normal bundle $\barN_{n, \beta} \rightarrow \barM_{0,n}(X, \beta)$,

\begin{theoremx}[Theorem \ref{thm:tropical_Euler_class}]
There exists a well-defined map that computes the tropical, inverse Euler class of the SUSY normal bundle $\barN_{n, \beta} \rightarrow \barM_{0,n}(X, \beta)$,

\begin{align*}
   e^{K, trop}(\barN_{n, \beta})^{-1}: A^{trop}_*(TSM_{0,n}(\Sigma_X, \beta)) &\rightarrow A^{trop}_*(TSM_{0,n}(\Sigma_X, \beta)) \otimes \mathbb{C}(\kappa)\\
   \alpha &\mapsto E^{trop}\left(\barN_{n, \beta}\right) \cap \alpha 
\end{align*}
where $E^{trop}\left(\barN_{n, \beta}\right)$ is defined in Definition \ref{eq:inverse_Euler_exp_trop}.
\end{theoremx} 

To show Theorem \ref{thm:tropical_Euler_class}, we make Assumption \ref{asu:locally_trivial}, which lays the groundwork for tropicalization of super Gromov-Witten invariants; it assumes that families of stable maps are in some sense locally trivialized by toroidal, spin varieties. Assumption \ref{asu:locally_trivial} allows us to define a procedure to compute the tropical, Euler class of the SUSY bundle $\barN_{n, \beta}$ concretely (see Remark \ref{rem:upshot of the assumption}). In Definition \ref{def:tropical_SGW}, we define tropical, genus-0, $n$ Neveu-Schwarz marked, super Gromov-Witten invariants of a convex, toric variety $X$. We use our procedure to compute a tropical, super Gromov-Witten invariant of $\bP^1$ in Example \ref{ex:SGW_P1}; namely, we show that the tropical, genus-0, 2-marked super Gromov-Witten invariant of $\bP^1$ in curve class $H \in H_2(\bP^1, \mathbb{Z})$ is $\frac{1}{\kappa^3}$.

We lastly show that our tropicalization of super Gromov-Witten invariants of a convex toric variety $X$ in Definition \ref{def:tropical_SGW} recovers Theorem \ref{thm:TSGW_point} when the target $X$ is a point,

\begin{corollaryx}[Corollary \ref{cor:connect}]
     Let $X = pt$. By evaluating the tropical, genus-0, $n$-marked, super Gromov-Witten invariant defined in Definition \ref{def:tropical_SGW} over the smooth locus $\left[\mathcal{M}_{0,n}\right] \subset \left[\barM_{0,n}\right]$, we recover Theorem \ref{thm:TSGW_point}.
\end{corollaryx}

\begin{remark} 
In \cite{KSY3}, the authors use localization to define and compute super Gromov-Witten invariants. In this paper, we work formally with the identity character $\kappa: \mathbb{C}^* \rightarrow \mathbb{C}^*$ to match expressions obtained by tropical geometry with localization calculations from \cite{KSY3}. 
\end{remark}

\subsection{Acknowledgments}
B.Z. thanks the Yau Mathematical Sciences Center, Tsinghua University and the Huiyan Talent Fund for support. A.S. is supported by grants from Beijing Institute of Mathematical Sciences and Applications (BIMSA), the Beijing NSF BJNSF-IS24005, and the China National Science Foundation (NSFC) NSFC-RFIS program W2432008. A.S. would like to also thank China's National Program of Overseas High Level Talent for generous support. 

\subsection{Notation}
\label{sec:notation}
Throughout this paper, $K := \mathbb{C}^*$ and $K' := \mathbb{C}^*$. Let $\kappa: K \rightarrow \mathbb{C}^*$ and $\kappa': K' \rightarrow \mathbb{C}^*$ be the identity characters. Let $\mathbb{C}(\kappa)$ be the field of fractions in $\kappa$ and similarly for $\kappa'$.

We write $\phi:\mathcal{C}/B \rightarrow X$ to  be a family of stable maps $\phi:\mathcal{C} \rightarrow X$, or the diagram

\begin{center}
    \begin{tikzcd}
    \mathcal{C} \arrow[d, "\pi"] \arrow[r, "\phi"] & X \\
    B
    \end{tikzcd}
\end{center}
Here $\pi:\mathcal{C} \rightarrow B$ is a family of pre-stable, genus-$g$, $n$-marked curves, i.e. a flat, projective map $\pi$ with $n$ sections $p_i: B \rightarrow \mathcal{C}$ such that for each $b \in B$, $\pi^{-1}(b)$ is a pre-stable, genus-$g$, $n$-marked curve $(C_b, p_1(b), \ldots, p_n(b))$.  For each $b \in B$, $\phi_b: (C_b, p_1(b), \ldots, p_n(b)) \rightarrow X$ is a stable map with $(\phi_b)_*[C_b] = \beta \in H_2(X, \mathbb{Z})$.

Define $rk_{n, \beta} := c_1(\beta) + n - 2$.

\section{SUSY normal bundles and super Gromov-Witten invariants}

\label{sec:SUSY_SGW}

Let $X$ be a convex, toric variety with $\dim X = k$. We recall the definition of the SUSY normal bundles $\barN_n \rightarrow \barM_{0,n}$ and $\barN_{n, \beta} \rightarrow \barM_{0,n}(X, \beta)$, and prove general facts about them in Section \ref{sec:SUSY_bundles}. We then review the definition of super Gromov-Witten invariants in Section \ref{sec:SGW}. 

\subsection{Spinor sheaf of nodal, genus-0 curve}
\label{sec:spinor sheaf}

We recall that for a smooth, genus-$g$, unmarked curve $C$, the spinor bundle $S_C$ is the line bundle satisfying $S_C^{\otimes 2} \cong K_C$. For example, when $C \cong \bP^1$, we have $S_C = \cO_{\bP^1}(-1)$. 

Suppose that $C$ is a genus-0, nodal curve with dual graph $\Gamma$. Let $\tilde{c}: \tilde{C} \rightarrow C$ be the normalization of $C$. The spinor sheaf $S_{\tilde{C}}$ of $\tilde{C}$ is the disjoint union of invertible sheaves $S_{\tilde{C}} := \bigsqcup_{i=1}^{V(\Gamma)} \cO_{\bP^1}(1)$ satisfying $S_{\tilde{C}}^{\otimes 2} = K^{\vee}_{\tilde{C}}$. We define the spinor sheaf $S_C$ of $C$ to be $S_C := \tilde{c}_* S_{\tilde{C}}$, which is rank 1 at smooth points and rank 2 at nodes. We do not take $S_C$ to be invertible in its definition.


\begin{remark}
   In this paper, we mainly work with the \textit{dual} spinor sheaf, i.e. a coherent sheaf $\mathcal{F}$ satisfying $\mathcal{F}^{\otimes 2} \cong K^{\vee}_X$ for a space $X$. 
\end{remark}

\subsection{SUSY normal bundles}

\label{sec:SUSY_bundles}
We recall the definition of the SUSY normal bundles $\barN_n \rightarrow \barM_{0,n}$ and $\barN_{n, \beta}\rightarrow \barM_{0,n}(X, \beta)$ from \cite{KSY3}. Let $\barM_{0,n}$ be the moduli space of genus-0, $n$-marked, stable curves \cite{HM}. For an $n$-pointed stable curve $(C, p_1,\ldots,p_n) \in \barM_{0,n}$, we denote its dual graph by $C_\Gamma$, or $\Gamma$ when there is no ambiguity. Let $\barM_{0,n}^{super}$ be the moduli space of genus-0, $n$ Neveu-Schwarz marked, stable supercurves \cite{KSY1} \cite{FKP}. Let $\barM^{super}_{0,n}(X, \beta)$ be the moduli space of genus-0, $n$ Neveu-Schwarz marked (Footnote \ref{fn:NS}), super stable maps constructed in \cite{KSY2}. From \cite{KSY2}, we have an embedding $\barM_{0,n}(X, \beta) \hookrightarrow \barM^{super}_{0,n}(X, \beta)$.

\subsubsection{SUSY normal bundle $\barN_n$}

Fix a dual graph $\Gamma$. Let $c: \mathcal{C} \rightarrow B$ be a family of genus-0, stable curves over $B$, where $\forall x \in B$, $\mathcal{C}_x := c^{-1}(x)$ is a stable curve with dual graph $\Gamma$. Let $b: B \rightarrow \barM_{0,n}$ be the classifying map of $c: \mathcal{C} \rightarrow B$. Assume that a dual spinor sheaf $S_{\mathcal{C}}$ of $\mathcal{C}$ exists, i.e. $S_{\mathcal{C}}^{\otimes 2} \cong K^{\vee}_{\mathcal{C}}$. For each half edge $t \in E(\Gamma)$, we have sections $p_t: B \rightarrow C$ such that $p_t(x)$ is a smooth marked point of $\mathcal{C}_x$, and for each compact edge $e \in E(\Gamma)$, we have sections $p_e: B \rightarrow C$ such that $p_e(x)$ is a node of $\mathcal{C}_x$.


Associated to $c: \mathcal{C} \rightarrow B$, let $\barN_{\mathcal{C}}$ be the coherent sheaf on $B$ defined by the short exact sequence,

\begin{equation}
\label{eq:SES_curves}
 0 \rightarrow c_*S_{\mathcal{C}} \rightarrow \bigoplus_{t \in E(\Gamma)}p_t^* S_\mathcal{C} \oplus \bigoplus_{e \in E(\Gamma)} p_e^* S_\mathcal{C} \rightarrow N_\mathcal{C} \rightarrow 0   
\end{equation}
The \textit{SUSY normal bundle} $\barN_{n}$ on $\barM_{0,n}$ is the coherent sheaf satisfying $b^* \overline{N}_{n} = N_{\mathcal{C}}$ for any family of curves $c: \mathcal{C} \rightarrow B$ with classifying map $b: B \rightarrow \barM_{0,n}$. 

\subsubsection{SUSY normal bundle $\barN_{n, \beta}$}

Fix a dual graph $\Gamma$. Let $\phi:\mathcal{C}/B \rightarrow X$ be a family of genus-0, stable maps over $B$, where for  $c: \mathcal{C}\rightarrow B$ and $\forall x \in B$, $\mathcal{C}_x := c^{-1}(x)$ is a pre-stable curve with dual graph $\Gamma$. Furthermore, $(\phi_x)_*[\mathcal{C}_x] = \beta \in H_2(X, \mathbb{Z})$ for all $x \in B$. Assume that a dual spinor sheaf $S_{\mathcal{C}}$ of $\mathcal{C}$ exists. The coherent sheaf $N_{(\mathcal{C}, \phi)}$ on $B$ is defined by the short exact sequence,

\begin{equation}
\label{eq:SES_maps}
0 \rightarrow c_* S_{\mathcal{C}} \rightarrow \bigoplus_{t \in E(\Gamma)}p_t^* S_{\mathcal{C}} \oplus \bigoplus_{e \in E(\Gamma)} p_e^* S_{\mathcal{C}} \oplus c_*(S_{\mathcal{C}}^{\vee} \otimes \phi^*TX) \rightarrow N_{(\mathcal{C}, \phi)} \rightarrow 0
\end{equation}
where the second map is defined by $s \mapsto \bigoplus_t p_t^*s \oplus \bigoplus_e p_e^* s  \oplus -\langle s, d\phi \rangle$. When $\phi$ is the product family of stable maps $\phi: (C \times U)/U \rightarrow X$ with $U = \mathbb{A}^n, \bP^n$, or a Calabi-Yau space, then the exact sequence \ref{eq:SES_maps} holds. When $X = pt$, by definition $\barN_{n, \beta} \cong \barN_n$. 

The \textit{SUSY normal bundle} $\barN_{n, \beta}$ on $\barM_{0,n}(X, \beta)$ is the coherent sheaf satisfying $b^* \overline{N}_{n, \beta} = N_{(C, \phi)}$ for any family of stable maps $\phi: \mathcal{C}/B \rightarrow X$ with classifying map $b: B \rightarrow \barM_{0,n}(X, \beta)$. The rank of $\barN_{n, \beta}$ is $rk_{n, \beta} := c_1(\beta) + n-2$. 

We show that under mild assumptions, the coherent sheaf $\barN_{n, \beta}$ is a vector bundle on $\barM_{0,n}(X, \beta)$.

\begin{proposition}
\label{prop:vector_bundle}
    Let $X$ be convex. Suppose we have a Kontsevich space $\barM_{0,n}(X, \beta)$ that admits a cover $\{B_i\}$ with classifying maps $b_i: B_i \rightarrow \barM_{0,n}(X, \beta)$ corresponding to families of stable maps $\phi_i: \mathcal{C}_i/B_i \rightarrow X$ such that $N_{(\mathcal{C}_i, \phi_i)}$ exists for all $i$. Then $\barN_{n, \beta}$ is a vector bundle on $\barM_{0,n}(X, \beta)$.

    \begin{proof}
        When $X$ is convex, $\barM_{0,n}(X, \beta)$ is a normal projective variety by \cite{FP}, Theorem 2, and hence it is reduced. Since $\barN_{n, \beta}$ is of constant rank over $\barM_{0,n}(X, \beta)$ by \cite{KSY3}, Lemma 2.3.5, it is locally free by \cite{Ha}, Exercise II.5.8.    
    \end{proof}
\end{proposition}

\begin{example}
    Let $B = pt$. Suppose we have $\phi: C/pt \rightarrow X$, where $C$ is a genus-0, $n$-pointed stable curve. Let $|V|$ be the number of vertices of $C$. By the Riemann-Roch formula, we have $H^0(S^{\vee}_C \otimes \phi^*TX) \cong \mathbb{C}^{c_1(\beta)}$. Hence,

    \[
    N_{(C, \phi)} \cong \left(\mathbb{C}^{n + 2(|V| - 1)} \oplus H^0(S_{C}^{\vee} \otimes \phi^* TX)\right)/\mathbb{C}^{2|V|} \cong \mathbb{C}^{n - 2} \oplus \mathbb{C}^{c_1(\beta)} \cong \mathbb{C}^{c_1(\beta) + n - 2}
    \]
\end{example}

\begin{example}
\label{ex:product_SUSY}
    We express the SUSY bundle for a product family of genus-0, $n$-marked stable maps $\phi_U:(\mathcal{C} \times U)/U \rightarrow X \in \barM_{0,n}(X, \beta)$ where $\mathcal{C} \subseteq \barM_{\Gamma}$ is a family of pre-stable curves of a fixed dual graph $\Gamma$, and the base $U$ is a spin variety. This means that $\forall u \in U$, $\mathcal{C}_u$ is a pre-stable, genus-0 curve of dual graph $\Gamma$, and we have $(\phi_{U, u})_*[\mathcal{C}_u] = \beta \in H_2(X, \mathbb{Z})$. By definition, there exists an invertible sheaf $S_U$ with $S_U^{\otimes 2} = K_U^{\vee}$. Let $0 \in U$ be the closed point.
    Denote the two projection maps by $\mathcal{C}_0 \xleftarrow[]{\pi} \mathcal{C} \times U  \xrightarrow[]{c} U$. 

    We recall from Section \ref{sec:spinor sheaf} that the spinor sheaf of $\mathcal{C}_0$ is defined as $S_{\mathcal{C}_0} := \tilde{c}_* S_{\tilde{\mathcal{C}}_0}$, where $S_{\tilde{\mathcal{C}}_0} = \bigsqcup_{i=1}^{|V(\Gamma)|} \cO_{\bP^1}(1)$. Since the fibers of $\mathcal{C}$ each have dual graph $\Gamma$, we have $S_{\mathcal{C}_i} \cong S_{\mathcal{C}_j}$ for $i, j \in U$. We define $S_{\mathcal{C} \times U} := \pi^* S_{C_0} \otimes c^* S_U$. Suppose that $\forall u_i, u_j \in U$, there exists an isomorphism $\varphi_{i,j}: \mathcal{C}_{u_i} \xrightarrow{\sim} \mathcal{C}_{u_j}$ such that $\phi_{U, u_i} = \phi_{U, u_j} \circ \varphi_{i,j}$. Then, we have $\phi_U^*TX \cong \pi^* \phi_{U, 0}^* TX$, where $\phi_{U, 0}: \mathcal{C}_0 \rightarrow X$ is the stable map from $\mathcal{C}_0$.
    
    Hence, from \ref{eq:SES_maps},  we have the isomorphism of sheaves,
    
    \begin{align*}
        N_{(\mathcal{C} \times U, \phi_U)} &\cong \left(\bigoplus_{t} p_t^{*} S_{\mathcal{C} \times U} \oplus \bigoplus_{e} p_e^* S_{\mathcal{C} \times U} \oplus c_*(S_{\mathcal{C} \times U}^{\vee} \otimes \phi_U^* TX)\right)/ c_*S_{\mathcal{C} \times U} \\
        &:= \left(\bigoplus_{t} p_t^{*} (\pi^* S_{\mathcal{C}_0} \otimes c^* S_U) \oplus \bigoplus_{e} p_e^* (\pi^* S_{\mathcal{C}_0} \otimes c^* S_U) \right. \\
        &\left. \oplus c_*((\pi^* S^{\vee}_{\mathcal{C}_0} \otimes c^* S^{\vee}_U) \otimes \pi^*\phi_{U, 0}^* TX)\right)/ c_*(\pi^*S_{\mathcal{C}_0} \otimes c^* S_U)
    \end{align*}
    Define $j_t := \pi \circ p_t$ and $j_e := \pi \circ p_e$. We have $c \circ p_t = Id$. By the projection formula, 

    \[
    N_{(\mathcal{C} \times U, \phi_U)} \cong \frac{\left(\bigoplus_{t} \left(j_t^* S_{\mathcal{C}_0} \otimes S_U \right) \oplus \bigoplus_{e} \left(j_e^* S_{\mathcal{C}_0} \otimes S_U\right) \oplus (c_*\pi^* (S^{\vee}_{\mathcal{C}_0} \otimes \phi_{U,0}^* TX) \otimes S_U^{\vee})\right)}{(c_* \pi^* S_{\mathcal{C}_0} \otimes S_U)}
    \]
    By definition, we have $j_e^* S_{\mathcal{C}_0} \cong \mathbb{C}^2, j_t^* S_{\mathcal{C}_0} \cong \mathbb{C}$, and $c_*\pi^* (S^{\vee}_{\mathcal{C}_0} \otimes \phi_{U,0}^* TX)$ is the constant sheaf on $U$ given by $\mathbb{C}^{c_1(\beta)}$. Hence, the above becomes,

    \[
     N_{(\mathcal{C} \times U, \phi_U)} \cong \left(S_U^{\oplus n} \oplus S_U^{\oplus 2(|V|-1)} \oplus (S_U^{\vee})^{\oplus c_1(\beta)} \right)/S_U^{\oplus 2|V|} \cong S_U^{\oplus n-2} \oplus (S_U^{\vee})^{\oplus c_1(\beta)}
    \]
    We have $c_1(N_{(\mathcal{C} \times U, \phi_U)}) = (n-2 - c_1(\beta))c_1(S_U)$. Hence, $N_{(\mathcal{C} \times U, \phi_U)}$ is expressed by the spinor bundle of $U$, and is independent of $\mathcal{C}$.
    
    For example, when $U = \bP^m$ with $m$ odd, and $S_U = \cO_{\bP^m}(\frac{m+1}{2})$ is the dual spinor bundle, then $N_{(\mathcal{C} \times U, \phi_{U})} \cong  \cO_{\bP^m}(\frac{m+1}{2})^{\oplus n-2} \oplus \cO_{\bP^m}(\frac{-m-1}{2})^{\oplus c_1(\beta)}$. If $U = \mathbb{A}^n$, $S_U = \cO_{\mathbb{A}^n}$, then $N_{(\mathcal{C} \times U, \phi_U)} \cong \cO_{\mathbb{A}^n}^{rk_{n, \beta}}$. If $U$ is Calabi-Yau, i.e. $K_U \cong \cO_U$, then we have $N_{(\mathcal{C} \times U, \phi_U)} \cong \cO_{U}^{rk_{n, \beta}}$ as well.
\end{example}

\subsection{Super Gromov-Witten invariants}
\label{sec:SGW}

We recall that a super Riemann surface $\Sigma$ is a complex $1\vert 1$-dimensional supermanifold $\Sigma$ with an odd holomorphic distribution $\mathcal{D} \subset T\Sigma$ such that $\frac{1}{2}[\cdot, \cdot]: \mathcal{D} \otimes_{\mathbb{C}} \mathcal{D} \cong T\Sigma/\mathcal{D}$ \cite{LR} \cite{Wit}. Endow $X$ with a symplectic form and compatible, almost complex structure $J$. In \cite{KSY1}, the authors define a genus-0, super $J$-holomorphic curve as a map $u: (\Sigma, \mathcal{D}, j) \rightarrow (X, J)$ from a genus-0, super Riemann surface $(\Sigma, \mathcal{D}, j)$ with almost complex structure $j$ to $(X, \omega, J)$, satisfying the equation,

\[
\frac{1}{2}(du + J \circ du \circ j)|_{\mathcal{D}} = 0
\] 
Super $J$-holomorphic curves, in analogy with $J$-holomorphic curves, are the critical points of a superconformal action. 

In \cite{KSY1}, the authors construct the moduli space $\barM^{super}_{0,n}(X, \beta)$ of genus-0, $n$ Neveu-Schwarz marked (Footnote \ref{fn:NS}), super stable maps to $X$ in curve class $\beta \in H_2(X, \mathbb{Z})$. In \cite{KSY2}, the authors define an action of $K' := \mathbb{C}^*$ on $\barM^{super}_{0,n}(X, \beta)$, with fixed locus the moduli space $\barM_{0,n}^{1/2}(X, \beta)$ of genus-0, $n$-marked 2-spin, stable maps studied in \cite{JKV1}. Let $\kappa': K' \rightarrow \mathbb{C}^*$ be the identity character of $K'$. Forgetting the 2-spin structure gives us an isomorphism of coarse moduli spaces $\barM_{0,n}^{1/2}(X, \beta) \cong \barM_{0, n}(X,\beta)$, which implies $H^*(\barM^{1/2}_{0,n}(X, \beta), \mathbb{Q}) \cong H^*(\barM_{0,n}(X, \beta), \mathbb{Q})$. Let $K := \mathbb{C}^*$ act trivially on $\barM_{0,n}(X, \beta)$, and $\kappa:K \rightarrow \mathbb{C}^*$ be the identity character of $K$. We have,

\begin{align*}
    H^*_{K'}(\barM_{0,n}^{1/2}(X,\beta), \mathbb{Q}) &\cong H^*(\barM_{0,n}^{1/2}(X,\beta), \mathbb{Q}) \otimes \mathbb{C}(\kappa') \\
    &\cong  H^*(\barM_{0,n}(X,\beta), \mathbb{Q}) \otimes \mathbb{C}(\kappa) \\
    &\cong H_K^*(\barM_{0,n}(X,\beta), \mathbb{Q}) 
\end{align*}

By Proposition \ref{prop:vector_bundle}, the SUSY normal bundle $\barN_{n, \beta}$ is a vector bundle on $\barM_{0,n}(X, \beta)$ of rank $c_1(\beta) + n - 2$. We let $K$ act on the fibers of $\barN_{n, \beta}$ by scaling. Denote $e^{K}(\barN_{n, \beta}) \in H^*_K(\barM_{0,n}(X, \beta), \mathbb{Q})$ to be the equivariant Euler class of $\barN_{n, \beta}$ with respect to $K$.

For $1 \leq i \leq n$, let $ev_i: \barM_{0,n}(X, \beta) \rightarrow X$ be the evaluation map at the $i$-th marking. Let $[\barM_{0,n}(X, \beta)]^{vir}$ be the virtual fundamental class of $\barM_{0,n}(X, \beta)$ \cite{BF}. In \cite{KSY3}, the genus-0, $n$ Neveu-Schwarz marked, super Gromov-Witten invariants of $X$ in class $\beta \in H_2(X, \mathbb{Z})$ with incidence conditions $\alpha_i \in H^*(X)$ are defined as,

\begin{definition}[\cite{KSY3}]
\label{def:SGW}
\begin{align*}
  SGW_{0,n}(X, \beta)(\alpha_1,\ldots, \alpha_n) = \int_{\left[\barM_{0,n}(X, \beta)\right]^{vir}}\frac{ev_1^*\alpha_1 \cup \ldots \cup ev_n^*\alpha_n}{e^K(\barN_{n,\beta})} \in H^*(\barM_{0,n}(X, &\beta), \mathbb{Q}) \\ &\otimes \mathbb{C}(\kappa)  
\end{align*}
where $e^K(\barN_{n, \beta})$ is the equivariant Euler class of $\barN_{n, \beta}$ with respect to $K := \mathbb{C}^*$. 
\end{definition}

\begin{remark}
    When $\sum_i \deg \alpha_i = \text{virdim }\barM_{0,n}(X, \beta)$, the super Gromov-Witten invariants are the ordinary Gromov-Witten invariants, modulo a power of $\kappa$ (see Equation \ref{eq:inverse_Euler_exp} for the expression of $e^K(\barN_{n, \beta})^{-1})$. 
    
    When $\sum_i \deg \alpha_i < \text{virdim }\barM_{0,n}(X, \beta)$, the ordinary Gromov-Witten invariant is trivial. However, there will be a potentially non-trivial contribution of $e^K(\barN_{n, \beta})^{-1}$, making the super Gromov-Witten invariant potentially non-trivial.
\end{remark}

\section{Tropicalization of super Gromov-Witten invariants of a point}

\label{sec:SGW_point}
In this section, we tropicalize super Gromov-Witten invariants of a point \cite{KSY3}. We first review aspects of the moduli of curves $\barM_{0,n}$ and its tropicalization $\barM^{trop}_{0,n}$. In Theorem \ref{thm:TSGW_point}, we show that the genus-0, super Gromov-Witten invariants of a point in Proposition \ref{prop:SGW_point_descendant} are equal to genus-0, tropical super Gromov-Witten invariants of a point defined in Definition \ref{def:TSGW(pt)}.

\subsection{Tropical geometry}

The tropical semi-ring $\mathbb{T}$ is the set $\mathbb{R} \cup \{\infty\}$ with two operations $\oplus$ and $\odot$, defined as $a\oplus b := \min(a,b)$ and $a \odot b := a+b.$ Let $\mathbb{T}^*$ be the multiplicative units of $\mathbb{T}$. Tropical projective space is defined as $\mathbb{T}\mathbb{P}^n := \{\lambda \in \mathbb{T}^n | \lambda \sim \lambda' \textit{ if }\exists \lambda'' \in \mathbb{T}^*, \lambda = \lambda' \odot \lambda''\}$. We have $\mathbb{T}\bP^1 = \overline{\mathbb{R}} = \mathbb{R} \cup \{\pm\infty\}$. 

\begin{definition}
\label{def:trop_curve}
A \textit{abstract tropical curve} is a graph $\Gamma$ and a weighted metric $w: E(\Gamma)\rightarrow \mathbb{R}_{\geq 0} \cup \{\infty\}$. We write $V(\Gamma)$ and $E(\Gamma)$ to be the vertex and edge sets of $\Gamma$ respectively, and $val(V)$ to be the valency of the vertex $V$.
\end{definition}

\subsection{Descendant invariants of $\barM_{0,n}$}

Let $\barM_{0,n}$ be the moduli space of genus-0, $n$-marked, stable curves \cite{HM}. For each $i$, we have the morphism $ft_i: \barM_{0,n+1} \rightarrow \barM_{0,n}$ which forgets the $i$-th marked point and stabilizes. For $1 \leq i \leq n$, denote $\mathbb{L}_i \rightarrow \barM_{0,n}$ to be the line bundle whose fiber above an $n$-pointed, stable curve $(C,p_1,\ldots, p_n)$ is $T^*_{p_i}C$. Let $\psi_i := c_1(\mathbb{L}_i) \in H^2(\barM_{0,n})$. Define $D_{i, n+1} \subseteq \barM_{0, n}$ to be the divisor consisting of rational curves with 2 components, with 1 component containing only the marked points $i, n+1$. We have,

\begin{lemma}[Comparison lemma, \cite{Koc}]
\label{lem:comparison_lemma}
For $1 \leq i \leq n$,

\[
\psi_i - ft^*_i \psi_i = D_{i, n+1}
\]
\end{lemma}
The $\psi_i$-classes satisfy $\psi_i \cdot D_{i, n+1} = \psi_{n+1} \cdot D_{i, n+1} = 0$. 

\subsection{Tropicalization of $\barM_{0,n}$}

We recall the construction of the moduli space $\barM_{0,n}^{trop}$ of genus-$0$, $n$-marked, tropical curves \cite{Abr} \cite{ACP} \cite{C}. The dual graph $\Gamma$ of a stable curve $C \in \barM_{0,n}$ is defined by introducing a vertex $V$ for each irreducible component of $C$, a half edge for each $p_i$, and a bounded edge connecting two vertices for each nodal point of $C$ (Figure \ref{def:trop_curve}). The dual graph $\Gamma$ of $C \in \barM_{0, n}$ is an abstract tropical curve in the sense of Definition \ref{def:trop_curve}. For each bounded edge $e \in E(\Gamma)$, introduce a cone $\mathbb{R}_{\geq 0}$ that tracks the length of $e$. For a graph $\Gamma$, denote $\barM_{\Gamma} \subset \barM_{0,n}^{trop}$ to be the subset of tropical curves with dual graph $\Gamma$. We see that $\barM_{\Gamma} \cong (\mathbb{R}_{\geq 0})^{|E(\Gamma)|}/Aut(\Gamma).$ We have an inclusion $\barM_{\Gamma'} \subset \barM_{\Gamma}$, if $\Gamma$ arises by contracting the length of a bounded edge in $\Gamma'$ to 0. The moduli space of tropical curves $\barM_{0,n}^{trop}$ is obtained by gluing along the inclusions $\barM_{\Gamma'} \subset \barM_{\Gamma}$, where $\Gamma$ is obtained by contracting $\Gamma'$. We have $\barM_{0,n}^{trop}$ is a cone complex of dimension $n-3$.

\begin{example}
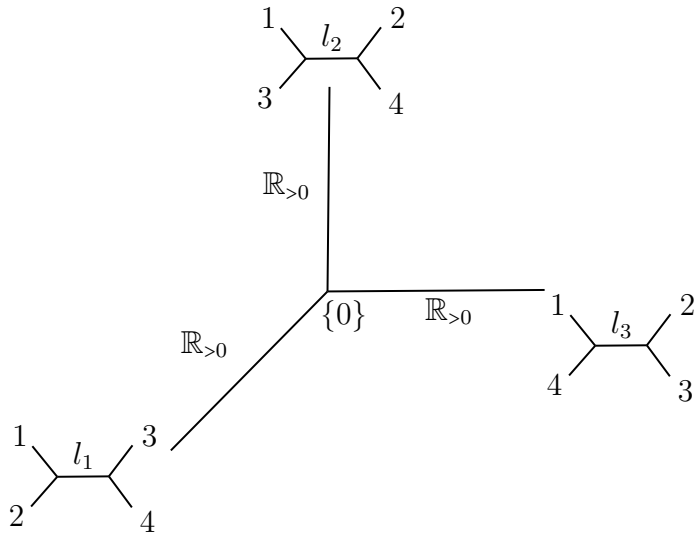
   The dual graphs of the 3 boundary divisors of $\barM_{0,4} \cong \bP^1$ contain a single compact edge whose length $l_i$ parametrizes a cone $\mathbb{R}_{\geq 0}$. Hence, $\barM_{0,4}^{trop}$ contains 3 copies of the cone $\mathbb{R}_{\geq 0}$ all glued at the $0$-cone.
\end{example}

\begin{definition}[\cite{KM} \cite{Mik}]
\label{def:tropical_psi}
The \textit{tropical Psi-class} $\Psi_i \subseteq \barM_{0,n}^{trop}$ is the fan of $(n-4)$-dimensional cones in $\barM_{0,n}^{trop}$ consisting of tropical curves whose half edge corresponding to the marked point $i$ is on a vertex $V$ with $val(V) = 4.$ 
\end{definition}

\begin{example}
    In Figure \ref{fig:super_dual_graph}, $C \in \Psi_1 \cap \Psi_2 \cap \Psi_4$.
\end{example}

\subsection{Super Gromov-Witten invariants of a point}
The super Gromov-Witten invariants of a point are given in Definition \ref{def:SGW} for $X = pt$. When $X = pt$, we have $\barN_{n, \beta} \cong \barN_n$. Hence, the genus-0, $n$-marked, super Gromov-Witten invariants of a point, $SGW_{0,n}(pt)$, are, 

\begin{definition}[\cite{KSY3}]
\label{def:SGW_point}
  \begin{align*}
        SGW_{0,n}(pt) := \int_{[\barM_{0,n}]}\frac{1}{e^K(\barN_n)} 
    \end{align*}
\end{definition}
For any $n \geq 4$ and $j \leq n$, there exists a line bundle $S_j \rightarrow \barM_{0,n}$ whose fiber above a curve $(C, p_1, \ldots, p_n) \in \barM_{0,n}$ is the vector space $V_j$ satisfying $V_j \otimes V_j = T_{p_j}C$. Let $ft_{n}:\barM_{0,n} \rightarrow \barM_{0,n-1}$ be the morphism that forgets the last marked point $p_n$ and stabilizes. From \cite{KSY3}, we have the splitting of $\barN_n \rightarrow \barM_{0,n}$ into line bundles,

\begin{equation}
\label{eq:Nn splitting}
    \barN_n \cong \mathbb{C} \oplus \bigoplus_{j=4}^n (ft^*)^{n-j}(S_j)
\end{equation}
where we have made an abuse of notation by writing $(ft^*)^{n-j}$ to mean the composition $ft^*_{n}\ldots ft^*_{j+1}$. Hence, from Equation \ref{eq:Nn splitting}, Definition \ref{def:SGW_point} becomes,

\[
SGW_{0,n}(pt) = \frac{(-1)^{n-3}}{2^{n-3}\kappa^{2n-5}}\sum_{\substack{i_4,\ldots, i_n \geq 0, \\ i_4 + \ldots + i_n = n-3}} \int_{\left[\barM_{0,n}\right]} (ft^*)^{n-4}\psi_4^{i_4} \ldots ((ft^*)\psi_{n-1}^{i_{n-1}})\psi_n^{i_n}
\]
where $\psi_l \in H^2(\barM_{0,l})$ for $4 \leq l \leq n$. 

Using Lemma \ref{lem:comparison_lemma}, we first re-express $SGW_{0,n}(pt)$ as descendant invariants on $\barM_{0,n}$,

\begin{proposition}
\label{prop:SGW_point_descendant}
    \begin{align*}
      SGW_{0,n}(pt) &= \frac{(-1)^{n-3}}{2^{n-3}\kappa^{2n-5}} \sum_{\substack{i_4,\ldots, i_n \geq 0, \\ i_4 + \ldots + i_n = n-3}}\int_{\barM_{0,n}} \psi_4^{i_4} \ldots \psi_{n-1}^{i_{n-1}}\psi_n^{i_n} \\
      &= \sum_{\substack{i_4,\ldots, i_n \geq 0, \\ i_4 + \ldots + i_n = n-3}} \frac{(n-3)!}{i_4!\ldots i_n!}  
    \end{align*}

    \begin{proof}
        We show that the integrand $(ft^*)^{n-4}\psi_4^{i_4} \ldots ((ft^*)\psi_{n-1}^{i_{n-1}})\psi_n^{i_n}$ is equal to,
        
        \[
        \psi_4^{i_4} \ldots \psi_{n-1}^{i_{n-1}}\psi_n^{i_n}
        \]
        where the $\psi_l$-classes for $4 \leq l \leq n$ are now understood to be in $H^2(\barM_{0,n})$.

        We first show,
        
        \[
        (ft^*(\psi_{n-1}^{i_{n-1}})) \psi_n^{i_n} =  \psi_{n-1}^{i_{n-1}} \psi_n^{i_n}
        \]
        By Lemma \ref{lem:comparison_lemma}, we have $ft^* \psi_{n-1} = \psi_{n-1} - D_{n-1, n}$. Since pullback commutes with cup product and $\psi_{n-1} \cdot D_{n-1, n} = 0$, then $ft^* (\psi_{n-1}^{i_{n-1}}) = (ft^* \psi_{n-1})^{i_{n-1}} = \psi_{n-1}^{i_{n-1}} + (-1)^{i_{n-1}}(D_{n-1, n})^{i_{n-1}}$. Since $\psi_{n} \cdot D_{n-1, n} = 0$, we have $ft^*(\psi_{n-1}^{i_{n-1}}) \psi_n^{i_n} =  \psi_{n-1}^{i_{n-1}} \psi_n^{i_n}$.

        By induction, suppose the following holds for $k \geq 1$,

        \[
        (ft^*)^k (\psi_{n-k}^{i_{n-k}})(ft^*)^{k-1} (\psi_{n-k+1}^{i_{n-k+1}}) \ldots \psi_{n}^{i_n} = \psi_{n-k}^{i_{n-k}}\psi_{n-k+1}^{i_{n-k+1}} \ldots \psi_{n}^{i_n}
        \]
        Consider,

        \[
        ((ft^*)^{k+1}(\psi_{n-k-1}^{i_{n-k-1}})) \psi_{n-k}^{i_{n-k}} \ldots \psi_n^{i_n}
        \]
        We have,

        \begin{align*}
            ((ft^*)^{k+1}(\psi_{n-k-1}^{i_{n-k-1}})) \psi_{n-k}^{i_{n-k}} \ldots \psi_n^{i_n} &= ((ft^*)^{k}(ft^*\psi_{n-k-1}^{i_{n-k-1}})) \psi_{n-k}^{i_{n-k}} \ldots \psi_n^{i_n} \\
            &= \left((ft^*)^k \left(\psi_{n-k-1}^{i_{n-k-1}} - D_{n-k-1, n-k}^{i_{n-k-1}}\right)\right)\psi_{n-k}^{i_{n-k}} \ldots \psi_n^{i_n} 
        \end{align*}
        where the second equality follows from Lemma \ref{lem:comparison_lemma}, and $\psi_{n-k-1} \cdot D_{n-k-1, n-k} = 0$. By the projection formula, $\left((ft^*)^k(D_{n-k-1, n-k}^{i_{n-k-1}})\right)\psi_{n-k}^{i_{n-k}} = D_{n-k-1, n-k}^{i_{n-k-1}} (ft_*)^k\psi_{n-k}^{i_{n-k}}$. However $(ft_*)^k: H^*(\barM_{0,n}) \rightarrow H^*(\barM_{0,n-k})$ will send $\psi_{n-k}$ to itself $\psi_{n-k}$. Since $\psi_{n-k} \cdot D_{n-k-1, n-k} = 0$, we have 
        
        \[
        \left((ft^*)^{k+1}(\psi_{n-k-1}^{i_{n-k-1}})\right) \psi_{n-k}^{i_{n-k}} \ldots \psi_n^{i_n} = \left((ft^*)^{k}(\psi_{n-k-1}^{i_{n-k-1}})\right) \psi_{n-k}^{i_{n-k}} \ldots \psi_n^{i_n}
        \]
        Hence, for all $k \geq 1$,
        
        \[
        \left((ft^*)^{k+1}(\psi_{n-k-1}^{i_{n-k-1}})\right) \psi_{n-k}^{i_{n-k}} \ldots \psi_n^{i_n} = \psi_{n-k-1}^{i_{n-k-1}} \psi_{n-k}^{i_{n-k}} \ldots \psi_n^{i_n}
        \]
        Setting $k = n-4$ proves the first desired equality.
        
        The second desired equality is a well-known consequence of Witten's Theorem, stated in \cite{HM}.
    \end{proof}
\end{proposition}

\subsubsection{Decorated dual graphs}

To define tropical, super Gromov-Witten invariants of a point, we first define a space $T_{0,n}$ of decorated dual graphs tracking powers of $\psi$-classes.
\begin{definition}
\label{def:super_dual_graphs}
    A \textit{decorated dual graph} $\Gamma_{\vec{k}} := (\Gamma, k_1,\ldots, k_n)$ consists of a dual graph $\Gamma \in \barM^{trop}_{0,n}$ of a genus-0, $n$-marked stable curve $C \in \barM_{0,n}$, and, for each half-edge of $\Gamma$, a non-negative integer $k_i \in \mathbb{Z}_{\geq 0}$ keeping track of of $\psi_i$. Define the space of decorated dual graphs $\overline{T}_{0,n} := \barM^{trop}_{0,n} \times \mathbb{Z}_{\geq 0}^n$. Let $T_{0,n} \subset \overline{T}_{0,n}$ be the subset of decorated dual graphs $\Gamma_{\vec{k}}$ where $\Gamma$ is the dual graph of a smooth curve.
\end{definition}

\begin{figure}[h!]
    \centering

\tikzset{every picture/.style={line width=0.75pt}} 

\begin{tikzpicture}[x=0.75pt,y=0.75pt,yscale=-1,xscale=1]

\draw   (100,136) .. controls (100,110.32) and (124.4,89.5) .. (154.5,89.5) .. controls (184.6,89.5) and (209,110.32) .. (209,136) .. controls (209,161.68) and (184.6,182.5) .. (154.5,182.5) .. controls (124.4,182.5) and (100,161.68) .. (100,136) -- cycle ;
\draw   (188,192.25) .. controls (188,171.13) and (203.67,154) .. (223,154) .. controls (242.33,154) and (258,171.13) .. (258,192.25) .. controls (258,213.37) and (242.33,230.5) .. (223,230.5) .. controls (203.67,230.5) and (188,213.37) .. (188,192.25) -- cycle ;
\draw    (284,151.03) -- (347,151.97) ;
\draw [shift={(349,152)}, rotate = 180.86] [color={rgb, 255:red, 0; green, 0; blue, 0 }  ][line width=0.75]    (10.93,-3.29) .. controls (6.95,-1.4) and (3.31,-0.3) .. (0,0) .. controls (3.31,0.3) and (6.95,1.4) .. (10.93,3.29)   ;
\draw [shift={(282,151)}, rotate = 0.86] [color={rgb, 255:red, 0; green, 0; blue, 0 }  ][line width=0.75]    (10.93,-3.29) .. controls (6.95,-1.4) and (3.31,-0.3) .. (0,0) .. controls (3.31,0.3) and (6.95,1.4) .. (10.93,3.29)   ;
\draw   (412,123) .. controls (412,116.92) and (417.6,112) .. (424.5,112) .. controls (431.4,112) and (437,116.92) .. (437,123) .. controls (437,129.08) and (431.4,134) .. (424.5,134) .. controls (417.6,134) and (412,129.08) .. (412,123) -- cycle ;
\draw    (384,81) -- (414,117) ;
\draw    (419,65.5) -- (424.5,112) ;
\draw    (369,111) -- (412,123) ;
\draw    (432.6,131.8) -- (475.6,178.8) ;
\draw   (472,186.2) .. controls (472,180.12) and (477.6,175.2) .. (484.5,175.2) .. controls (491.4,175.2) and (497,180.12) .. (497,186.2) .. controls (497,192.28) and (491.4,197.2) .. (484.5,197.2) .. controls (477.6,197.2) and (472,192.28) .. (472,186.2) -- cycle ;
\draw    (484.5,197.2) -- (483.4,233.4) ;
\draw    (495,192.8) -- (526.6,224.2) ;
\draw    (128,136) ;
\draw [shift={(128,136)}, rotate = 0] [color={rgb, 255:red, 0; green, 0; blue, 0 }  ][fill={rgb, 255:red, 0; green, 0; blue, 0 }  ][line width=0.75]      (0, 0) circle [x radius= 3.35, y radius= 3.35]   ;
\draw    (174.6,139.8) ;
\draw [shift={(174.6,139.8)}, rotate = 0] [color={rgb, 255:red, 0; green, 0; blue, 0 }  ][fill={rgb, 255:red, 0; green, 0; blue, 0 }  ][line width=0.75]      (0, 0) circle [x radius= 3.35, y radius= 3.35]   ;
\draw    (152.5,163.75) ;
\draw [shift={(152.5,163.75)}, rotate = 0] [color={rgb, 255:red, 0; green, 0; blue, 0 }  ][fill={rgb, 255:red, 0; green, 0; blue, 0 }  ][line width=0.75]      (0, 0) circle [x radius= 3.35, y radius= 3.35]   ;
\draw    (214.5,195.25) ;
\draw [shift={(214.5,195.25)}, rotate = 0] [color={rgb, 255:red, 0; green, 0; blue, 0 }  ][fill={rgb, 255:red, 0; green, 0; blue, 0 }  ][line width=0.75]      (0, 0) circle [x radius= 3.35, y radius= 3.35]   ;
\draw    (240.4,208.25) ;
\draw [shift={(240.4,208.25)}, rotate = 0] [color={rgb, 255:red, 0; green, 0; blue, 0 }  ][fill={rgb, 255:red, 0; green, 0; blue, 0 }  ][line width=0.75]      (0, 0) circle [x radius= 3.35, y radius= 3.35]   ;

\draw (122,110.4) node [anchor=north west][inner sep=0.75pt]    {$p_{1}$};
\draw (142,137.9) node [anchor=north west][inner sep=0.75pt]    {$p_{2}$};
\draw (167,117.4) node [anchor=north west][inner sep=0.75pt]    {$p_{4}$};
\draw (203,171.4) node [anchor=north west][inner sep=0.75pt]    {$p_{3}$};
\draw (233,185.4) node [anchor=north west][inner sep=0.75pt]    {$p_{5}$};
\draw (325.5,115.9) node [anchor=north west][inner sep=0.75pt]    {$( p_{1} ,\ k_{1})$};
\draw (407.1,36.9) node [anchor=north west][inner sep=0.75pt]    {$( p_{4} ,\ k_{4})$};
\draw (332,56.4) node [anchor=north west][inner sep=0.75pt]    {$( p_{2} ,\ k_{2})$};
\draw (515.1,189.9) node [anchor=north west][inner sep=0.75pt]    {$( p_{5} ,\ k_{5})$};
\draw (416.6,221.4) node [anchor=north west][inner sep=0.75pt]    {$( p_{3} ,\ k_{3})$};
\draw (150.67,51.07) node [anchor=north west][inner sep=0.75pt]    {$C$};

\end{tikzpicture}

    \caption{Decorated dual graph for $C \in \barM_{0,5}$.}
    \label{fig:super_dual_graph}
\end{figure}
Let $A^{trop}_*(\barM^{trop}_{0,n})$ be the Chow group of tropical cycles of $\barM^{trop}_{0,n}$ \cite{CGM}. Let $[\barM_{0,n}^{trop}]$ be the top-dimensional tropical cycle on $\barM_{0,n}^{trop}$ that is the constant function 1  \cite{CG}, Definition 5.9.

\begin{definition}
\label{def:W}
    Define the function, 

    \begin{align*}
         W: T_{0,n} &\rightarrow A^{trop}_*(\barM^{trop}_{0,n}) \otimes \mathbb{C}(\kappa) \\
         \Gamma_{\vec{k}} = (\Gamma, k_1,\ldots, k_n) &\mapsto \frac{(-1)^{n-3}}{2^{n-3}\kappa^{2n-5}}\int_{[\barM^{trop}_{0,n}]} \Psi_4^{k_4}\ldots\Psi_n^{k_n}
    \end{align*}
    The integral on $\barM^{trop}_{0,n}$ is computed via tropical intersection \cite{CGM}.
\end{definition}
Let $Z := \{\Gamma_{\vec{k}} = (\Gamma, k_1,\ldots, k_n) | k_4 + \ldots k_n = n-3\} \subset T_{0,n}$, where $\Gamma$ is the dual graph of a smooth curve.

\begin{definition}
\label{def:TSGW(pt)}
    The tropical, genus-0, $n$-marked super Gromov-Witten invariant of a point, $TSGW_{0,n}(pt)$, is,

    \[
    TSGW_{0,n}(pt) := \sum_{\Gamma_{\vec{k}} \in Z}W(\Gamma_{\vec{k}})
    \]
    where $W$ is defined in Definition \ref{def:W}. 
\end{definition}

\begin{theorem}
\label{thm:TSGW_point}
    Let $X = pt$. We have,
 
    \[
    SGW_{0,n}(pt) = TSGW_{0,n}(pt)
    \]
    where $SGW_{0,n}(pt)$ is given in Proposition \ref{prop:SGW_point_descendant}, and $TSGW_{0,n}(pt)$ is defined in Definition \ref{def:tropical_SGW}.
    
    \begin{proof}
    By \cite{KM}, Theorem 4.1, the tropical curves lying in the locus $\Psi_4^{k_4}\ldots\Psi_n^{k_n}$ for $k_4 + \ldots + k_n = n-3$ are the dual graphs of smooth curves. By \cite{KM}, Corollary 4.2, we have $\Psi_4^{k_4}\ldots\Psi_n^{k_n} = \psi_4^{k_4}\ldots\psi_n^{k_n}$, with $k_4 + \ldots + k_n = n-3$.  By Proposition \ref{prop:SGW_point_descendant} and Definition \ref{def:TSGW(pt)}, the desired claim follows. 
    \end{proof}
\end{theorem}

\section{Tropicalization of super Gromov-Witten invariants of a convex, toric variety} 

\label{sec:SGW_toric}

Let $X$ be a convex, toric variety with $\dim X = k$. In this section, we define and compute tropical, super Gromov-Witten invariants of $X$. In Section \ref{sec:tropical_lb}, we first review tropical Chern classes and tropicalization \cite{CG} \cite{CGM}. In Section \ref{sec:def_tropical_Euler}, we make an Assumption \ref{asu:locally_trivial}, stating that families of stable maps are in some sense locally trivialized by toroidal, spin varieties. Using Assumption \ref{asu:locally_trivial}, we define the tropical, inverse Euler class of $\barN_{n, \beta}$. Then, in Definition \ref{def:tropical_SGW}, we define tropical, genus-0, $n$-marked, super Gromov-Witten invariants of $X$. In Example \ref{ex:SGW_P1}, we compute a tropical, super Gromov-Witten invariant of $\bP^1$. In Corollary \ref{cor:connect}, we show our definition of tropical, genus-0, super Gromov-Witten invariants of $X$ in Definition \ref{def:tropical_SGW} recovers Theorem \ref{thm:TSGW_point} when $X = pt$. 

\subsection{Tropical line bundles}

\label{sec:tropical_lb}
A tropical space is a topological space $B$ with a sheaf $Aff_B$ of $\overline{\mathbb{R}}$-valued, piecewise-linear, integral affine functions. Examples of tropical spaces include $\mathbb{T}^n$ and $\mathbb{T}\bP^n$. From \cite{CGM}, Section 6.1.3, a tropical line bundle $L$ over a tropical space $B$ has 3 equivalent descriptions: 1) locally $L$ is given by $U_i \times \mathbb{T} \rightarrow U_i$ for an open cover $\{U_i\}$ of $B$ 2) an $Aff_B$-torsor 3) a cocycle $\xi \in H^1(B, Aff_B)$. Tropical line bundles $\cO_{\mathbb{T}\bP^1}(m) \rightarrow \mathbb{T}\bP^1$ are described in \cite{CGM}, Example 6.8.

\subsection{Tropicalization of line bundles}
\label{sec:tropicalization_lb}

Given a toroidal variety $X$ with no self intersections, tropicalization \cite{CG} produces a cone complex $\Sigma_X$ from $X$, which is a tropical space. For example, if $X$ is a toric variety, then $\Sigma_X$ is the fan of $X$. From a toric variety $X$, we denote by $\overline{\Sigma}_X$ the corresponding tropical toric variety \cite{Kaj} \cite{Pay}. For example, if $X = \bP^2$, then $\overline{\Sigma}_X = \mathbb{T}\bP^2$. A morphism of toroidal varieties $f: X \rightarrow Y$ induces a morphism of cone complexes $Trop(f): \Sigma_X \rightarrow \Sigma_Y$ induced by pullback by $f^*$ of Cartier divisors \cite{Gro2}. If we have a Chow cycle $c \in A^*(X)$, then the tropicalization of $c$ is the function $Trop(c): \Sigma_X \rightarrow \mathbb{Z}$ defined by $\sigma \mapsto \int_X c \cdot V(\sigma)$, where $V(\sigma)$ is the strata in $X$ corresponding to  a cone $\sigma \subset \Sigma_X$.

Let $L \rightarrow X$ be a line bundle. In \cite{CG}, the authors define tropicalization of $L$, which produces a tropical line bundle on the cone complex $\Sigma_X$. We say $L$ is \textit{tropicalizable} if its tropicalization is a tropical line bundle. We denote tropicalizations of a toroidal variety $X$ and of a line bundle $L \rightarrow X$ as $Trop(X)$ and $Trop(L)$, respectively. 

\begin{example}
    Let $X$ be the toroidal variety $\barM_{0,n}$. The tropicalization of $X$ is $\barM_{0,n}^{trop}$ by \cite{CG}, Theorem C.
\end{example} 

\begin{example}
\label{ex:spinor_bundle_g=0}
   Consider the toroidal variety $\bP^1$ with the anticanonical divisor $2[pt]$. We have $\Sigma_{(\bP^1, 2H)} \cong \mathbb{T}\bP^1$, and $Trop(\cO_{\bP^1}(m)) = \cO_{\mathbb{T}\bP^1}(m)$ by \cite{CG}, Example 4.3.
\end{example}

\subsection{Tropicalization of spinor sheaf}

\label{sec:tropical_spinor_sheaf}
Let $C$ be a genus-0, nodal curve with dual graph $\Gamma$, and let $\tilde{c}: \tilde{C} \rightarrow C$ be the normalization.
The spinor sheaf of the normalization $\tilde{C}$ is $S_{\tilde{C}} := \bigsqcup_{i=1}^{|V(\Gamma)|} \cO_{\bP^1}(1)$. Define $Trop(S_{\tilde{C}}) := \bigsqcup_{i=1}
^{|V(\Gamma)|} Trop(\cO_{\bP^1}(1))$. Write an element $\varphi \in Trop(S_{\tilde{C}})$ as $\varphi = (\varphi_1, \ldots \varphi_{|V(\Gamma)|})$. Let $K := \{ \varphi \in Trop(S_{\tilde{C}}) | \varphi(p_i) = \varphi(q_i), \forall i\} \subseteq Trop(S_{\tilde{C}})$. Let $Aff_{\tilde{C}} := \bigsqcup_{i=1}^{|V(\Gamma)|} Aff_{\bP^1}$. Define the tropical spinor sheaf of a genus-0, nodal curve $C$ as $Trop(S_C) := \tilde{c}_* K$.

\begin{proposition}
\label{prop:trop_of_S_C}
    $Trop(S_C)$ is an $Aff_{\tilde{C}}$-torsor.

    \begin{proof}
        First, we see that in neighborhoods of smooth points of $C$, $Trop(S_C)$ is an $Aff_{\bP^1}$-torsor by Example \ref{ex:spinor_bundle_g=0}, and hence it is an $Aff_{\tilde{C}}$-torsor. Let $U_i$ be a neighborhood of a node $n_i \in C$ with preimages $\tilde{c}^{-1}(n_i) = \{p_i, q_i\}$. We check that $\Gamma(U_i, Trop(S_C))$ is an $Aff_{\tilde{C}}$-torsor.  We have $\Gamma(U_i, Trop(S_C)) = \Gamma(U_{i1} \sqcup U_{i2}, K)$, where $U_{i1}, U_{i2} \subset \tilde{C}$ are neighborhoods of $p_i, q_i$ respectively. Suppose we have an affine function $s_{i1} \sqcup s_{i2}$ on $U_{i1} \sqcup U_{i2}$ such that $s_{i1}(p_i) = s_{i2}(q_i)$. Given an affine function $(\ell_{i1}, \ell_{i2})$ on $U_{i1} \sqcup U_{i2}$ such that $\ell_{i1}(p) = \ell_{i2}(q)$, there exists a unique $v_{i1} \sqcup v_{i2} \in Aff(\tilde{C})$ such that $s_{ij} + v_{ij} = \ell_{ij}$ for $j = 1,2$, since $Trop(\cO_{\bP^1}(1))$ is by definition an $Aff_{\bP^1}$-torsor \cite{CGM}.
        \end{proof}
\end{proposition}

\subsection{Tropical Chern classes}

\label{sec:tropical_Chern}
For tropical spaces $B$, Chow groups $A^{trop}_*(B)$ of tropical cycles are defined in \cite{CGM}, Definition 6.6. The tropical, first Chern class is a map $A^{trop}_*(B) \rightarrow A^{trop}_{*-1}(B)$ defined by intersection with a tropical Cartier divisor. For details on tropical intersection theory, we refer to \cite{CGM}, Section 6.1.2., and \cite{AR}. 

\begin{example}
\label{ex:chern_class_O(m)}

Let $X$ be the toroidal variety given $\bP^1$ and its anticanonical divisor $2[pt]$. The extended cone complex $\overline{\Sigma}_X$ of $X$ consists of the extended cones $\{0\}, \sigma_1 := \mathbb{R}_{\geq 0}(1,0) \cup \{\infty\}, \sigma_2 := \mathbb{R}_{\geq 0}(-1,0) \cup \{-\infty\}$. We have $\overline{\Sigma}_X \cong \mathbb{T}\bP^1 \cong \mathbb{R} \cup \{\pm \infty\}$.  Consider the line bundle $\cO_{\bP^1}(m) \rightarrow \bP^1$. Using tropicalization, we have,

\begin{center}
\begin{tikzcd}
\label{diag:trop_Chern}
\cO_{\bP^1}(m)  \arrow[r, "Trop"]\arrow[d] & Trop(\cO_{\bP^1}(m)) \arrow[d] \\
X  \arrow[r, "Trop"] & \overline{\Sigma}_X = [-\infty, \infty]
\end{tikzcd}
\end{center}
If $U_0 = [-\infty, \infty), U_1 = (-\infty, \infty] \subset \overline{\Sigma}_X$, then $Trop(\cO_{\bP^1}(m))$ is the tropical line bundle on $\overline{\Sigma}_X$ with local data $(U_0, mx)$ and $(U_{1}, 0)$. We see that $Trop(\cO_{\bP^1}(m))$ is isomorphic to $\cO_{\mathbb{T}\bP^1}(m)$ from \cite{CGM}, Example 6.8. The tropical first Chern class of $\cO_{\mathbb{T}\mathbb{P}^1}(m)$ is a map of tropical cycles $A^{trop}_*(\mathbb{T}\bP^1) \rightarrow A^{trop}_{*-1}(\mathbb{T}\bP^1)$ defined by $\alpha \mapsto D \cdot \alpha$, where $D = div(\phi)$ is the tropical Cartier divisor \cite{CGM} defining $\cO_{\mathbb{T}\mathbb{P}^1}(m)$. Consider the tropical fundamental class $[\overline{\Sigma}_X]$, which is the tropical 1-cycle given by the constant function 1 on $\overline{\Sigma}_X$ (\cite{CG}, Definition 5.9). The result of  $c_1(\cO_{\mathbb{T}\bP^1}(m))\cdot [\overline{\Sigma}_X]$ is the function on $\overline{\Sigma}_X$ with value $m$ at the origin and 0 elsewhere.


\end{example}

\subsection{Tropical, inverse Euler class of the SUSY bundle}

\label{sec:def_tropical_Euler}

In this section, we define the tropical, inverse Euler class $e^{K, trop}(\barN_{n, \beta})^{-1}$ of the SUSY bundle $\barN_{n, \beta}$ using techniques from \cite{CGM} \cite{CG}. In Definition \ref{def:tropical_Euler_local}, we first define the tropical, inverse Euler class for a product family of stable maps $\phi_{U_i}: (\mathcal{C}_{U_i} \times U_i)/U_i \rightarrow X$ satisfying Assumption \ref{asu:locally_trivial}. Then, assuming we can cover the Kontsevich space $\barM_{0,n}(X, \beta)$ by a collection $\{\phi_{U_i}:(\mathcal{C}_{U_i} \times U_i)/U_i \rightarrow X\}_i$ each satisfying Assumption \ref{asu:locally_trivial}, we extend the definition of tropical, inverse Euler class to dual graph strata $\barM_{\Gamma}(X, \beta) \subseteq \barM_{0,n}(X, \beta)$ in Lemma \ref{lem:strata}, and then to $\barM_{0,n}(X, \beta)$ in Theorem \ref{thm:tropical_Euler_class}.

Let $TSM_{0,n}(\Sigma_X, \beta)$ be the extended cone complex of genus-0, tropical stable maps \cite{Ran}. We make the following assumption about families of stable maps to define the tropical, inverse Euler class of $\barN_{n, \beta}$,

\begin{assumption}
\label{asu:locally_trivial}
    Let $\phi: \mathcal{C}/B \rightarrow X$ be a family of genus-0, stable maps in curve class $\beta$ over a base space $B$, where the domain curves $\mathcal{C}_b$ are of dual graph $\Gamma$ for all $b \in B$. Suppose $B$ admits a covering $\{U_i\}$ such that each $U_i$ is a toroidal, spin variety. Let $S_{U_i}$ be the invertible sheaf on $U_i$ such that $S_{U_i}^{\otimes 2} = K^{\vee}_{U_i}$. Assume that $\phi: \mathcal{C}/B \rightarrow X$ locally trivializes on each $U_i$, in the sense that $\phi_{U_i}$ is the product family $(\mathcal{C}_{U_i} \times U_i)/U_i \rightarrow X$, where $\mathcal{C}_{U_i} \rightarrow U_i$ is a family of pre-stable, genus-0 curves of dual graph $\Gamma$, and $S_{U_i}$ is a tropicalizable line bundle on $U_i$, i.e. $Trop(S_{U_i})$ is a tropical line bundle on the cone complex $\Sigma_{U_i}$ associated to $U_i$ \cite{CG} \cite{CGM}. Assume that the defining exact sequence \ref{eq:SES_maps} for $N_{(\mathcal{C}_{U_i} \times U_i, \phi_{U_i})}$ is defined on each $U_i$. Let $b_i:U_i \rightarrow \barM_{0,n}(X, \beta)$ be the classifying maps of $\phi_{U_i}$. Assume that for all $i$, $Trop(b_i):\Sigma_{U_i} \rightarrow TSM_{0,n}(\Sigma_X, \beta)$ is proper, and pulls back tropical cycles to tropical cycles, i.e. $Trop(b_i)^*: A_*^{trop}(TSM_{0,n}\Sigma_X, \beta)) \rightarrow A_*^{trop}(\Sigma_{U_i})$ is a well-defined map.
\end{assumption}
In Example \ref{ex:SGW_P1}, we have an example of a Kontsevich space $\barM_{0,2}(\bP^1, H)$ whose dual graph strata of stable maps satisfy Assumption \ref{asu:locally_trivial}.


\begin{example}
    We give an example of a family of stable maps $\phi:\mathcal{C}/B \rightarrow X$ satisfying Assumption \ref{asu:locally_trivial}. Let $B = \bP^2$ with its standard toric structure. Take the torus orbits $U_0 := (\mathbb{C}^*)^2$ and $U_i = \bP^1$ for $i = 1,2,3$. We see that $\{U_i\}_i$ is a stratification of $\bP^2$ by toroidal, spin varieties. Let $\phi:\mathcal{C}/\bP^2 \rightarrow X$ be a family of stable maps in curve class $\beta$ that trivializes over each $U_i$ in the sense of Assumption \ref{asu:locally_trivial}. This means we have product families of stable maps $\phi_{U_i}: (\mathcal{C}_{U_i} \times U_i)/U_i \rightarrow X$ for all $i$, and $N_{\left(\mathcal{C}_{U_i} \times U_i, \phi_{U_i}\right)}$ is defined on all $U_i$. By Example \ref{ex:product_SUSY}, we have $N_{(\mathcal{C}_{U_0} \times U_0, \phi_{U_0})} \cong S_{(\mathbb{C}^*)^2}^{\oplus c_1(\beta) + n-2}  \cong \cO_{(\mathbb{C}^*)^2}^{\oplus c_1(\beta) + n-2}$. On $U_i \cong \bP^1$, we have $N_{(\mathcal{C}_{U_i} \times U_i, \phi_{U_i})} \cong \cO_{\bP^1}(1)^{\oplus n-2} \oplus \cO_{\bP^1}(-1)^{\oplus c_1(\beta)}$ . Hence $N_{(\mathcal{C}_{U_i} \times U_i, \phi_{U_i})}$ is composed of tropicalizable line bundles for all $U_i$ by Example \ref{ex:spinor_bundle_g=0}. Since $U_i \cap U_j = pt$ for $i, j \geq 1, i \neq j$, we have $N_{(\mathcal{C}_{U_i} \times U_i, \phi_{U_i})}|_{U_i \cap U_j} \cong N_{(\mathcal{C}_{U_j} \times U_j, \phi_{U_j})}|_{U_i \cap U_j} \cong \mathbb{C}^{\oplus c_1(\beta) + n-2}$.
\end{example}

\subsubsection{The inverse equivariant Euler class}
Recall from Proposition \ref{prop:vector_bundle} that $\barN_{n, \beta} \rightarrow \barM_{0,n}(X, \beta)$ is a vector bundle of rank $rk_{n, \beta}$. Let $K := \mathbb{C}^*$ act trivially on the base $\barM_{0,n}(X, \beta)$, and by scaling on the fibers of $\barN_{n, \beta}$. We compute the inverse, equivariant Euler class $e^K(\barN_{n, \beta})^{-1}$ of $\barN_{n, \beta}$ with respect to $K$. Let $\kappa:K \rightarrow \mathbb{C}^*$ be the identity character. By the splitting principle, we write $\barN_{n, \beta} = \bigoplus_{i=1}^{rk_{n, \beta}} N_i$, where the $N_i$ are line bundles. For a vector bundle $V$ that splits $V = \bigoplus_{i=1}^m L_i$ into a direct sum of $m$ line bundles $L_i$, define,

\begin{equation}
\label{eq:inverse_Euler_exp}
    E\left(V\right) := \sum_{j \geq 0}\frac{(-1)^j}{\kappa^{m+j}}\sum_{\substack{i_1 \geq 0, \ldots, i_{m} \geq 0, \\ i_1 + \ldots + i_{m} = j}}c_1(L_1)^{i_1} \cup \ldots \cup c_1(L_m)^{i_{m}} 
\end{equation}
By \cite{KSY3}, Equation 3.2.1., 

\begin{equation}
\label{eq:inverse_euler}
 e^K(\barN_{n, \beta})^{-1} = E\left(\barN_{n, \beta}\right) 
\end{equation}
We consider $c_1(L_i)$ as a map $A_*(\barM_{0,n}(X, \beta), \mathbb{Q}) \rightarrow A_{*-1}(\barM_{0,n}(X, \beta), \mathbb{Q})$. When the $L_i$ are tropicalizable, define,

\begin{equation}
\label{eq:inverse_Euler_exp_trop}
    E^{trop}\left(V\right) := \sum_{j \geq 0}\frac{(-1)^j}{\kappa^{m+j}}\sum_{\substack{i_1 \geq 0, \ldots, i_{m} \geq 0, \\ i_1 + \ldots + i_{m} = j}}c_1(Trop(L_1))^{i_1} \cup \ldots \cup c_1(Trop(L_m))^{i_{m}} 
\end{equation} 
where the tropical first Chern class is taken in Equation \ref{eq:inverse_Euler_exp_trop}. 

\subsubsection{Tropical, inverse Euler class for the product family}
Suppose $\phi_U:(\mathcal{C}_U \times U)/U \rightarrow X$ is a product family satisfying Assumption \ref{asu:locally_trivial}, i.e. $U$ is a toroidal, spin variety, and let $N_{(\mathcal{C}_U \times U, \phi_U)}$ be the SUSY normal bundle associated to $\phi_U$. Let $b_U: U \rightarrow \barM_{0,n}(X, \beta)$ the classifying map of $\phi_U$. By definition, $N_{(\mathcal{C}_U \times U, \phi_U)} = b_U^* \barN_{n, \beta}$. We have the splitting $N_{(\mathcal{C}_U \times U, \phi_U)} = \bigoplus_{i=1}^{rk_{n, \beta}} N'_i$, where $N_i' := b_U^* (N_i)$. 

\begin{definition}
\label{def:tropical_Euler_local}
The tropical, inverse Euler class $e^{K, trop}(\barN_{(\mathcal{C}_U \times U, \phi_U)})^{-1}$ of the SUSY normal bundle $\barN_{(\mathcal{C}_U \times U, \phi_U)} = \bigoplus_{i=1}^{rk_{n, \beta}} N'_i$ for a product family  $\phi_U:(\mathcal{C}_U \times U)/U \rightarrow X$ satisfying Assumption \ref{asu:locally_trivial} is,

\begin{align*}
   e^{K, trop}(\barN_{(\mathcal{C}_U \times U, \phi_U)})^{-1}: A^{trop}_*(\Sigma_U) &\rightarrow A^{trop}_*(\Sigma_U) \otimes \mathbb{C}(\kappa)\\
   \alpha &\mapsto E^{trop}\left(\barN_{(\mathcal{C}_U \times U, \phi_U)}\right) \cap \alpha
\end{align*}
where the RHS is defined by intersection of tropical cycles \cite{CGM} \cite{CG}, and $E^{trop}(\cdot)$ is defined in Equation \ref{eq:inverse_Euler_exp_trop}.
\end{definition}

\begin{example}
\label{ex:trop_Euler_class}
    Suppose $U$ is the toroidal, spin variety $U = (\bP^1, 2H)$ with dual spinor bundle $S_{\bP^1} = \cO_{\bP^1}(1)$. We compute the tropical Euler class of $N_{(\mathcal{C}_U \times U, \phi_{U})}$ associated to the product family  $\phi_{U}: (\mathcal{C}_U \times U)/U \rightarrow X$. From Example \ref{ex:product_SUSY}, we have the splitting,
    
    \begin{equation}
    \label{eq:splitting}
       N_{(\mathcal{C}_U \times U, \phi_{U})} \cong \cO_{\bP^1}(1)^{\oplus n-2} \oplus \cO_{\bP^1}(-1)^{\oplus c_1(\beta)} 
    \end{equation}
    By definition, $N_{(\mathcal{C}_U \times U, \phi_{U})} \cong b_{U}^* \barN_{n, \beta}$. Hence, $b_{\bP^1}^* N_i$ is isomorphic to either $\cO_{\bP^1}(1)$ or $\cO_{\bP^1}(-1)$. By Example \ref{ex:spinor_bundle_g=0}, $Trop(\cO_{\bP^1})(\pm 1)$ are tropical line bundles on $\mathbb{T}\bP^1 = \overline{\Sigma}_{(\bP^1, 2H)}$. The tropical Chern class of $Trop(\cO_{\bP^1}(\pm 1))$ with respect to the toroidal variety $(\bP^1, 2H)$ is given in Example \ref{ex:chern_class_O(m)}. Hence, the tropical Euler class of $N_{(C \times \bP^1, \phi_{\bP^1})}$ is also obtained.
\end{example}

\subsubsection{Tropical, inverse Euler class of $\barN_{n, \beta}$}
\label{sec:tropical_Euler_class}
We recall that $\barM_{0,n}(X, \beta)$ is stratified by dual graph, i.e. $\barM_{0,n}(X, \beta) = \bigcup_{\Gamma} \barM_{\Gamma}(X , \beta)$, where we denote $\barM_{\Gamma}(X , \beta)$ to be the sub-moduli space of stable maps with domain curve of dual graph $\Gamma$. Recall that $TSM_{0,n}(\Sigma_X, \beta)$ is the extended cone complex of genus-0, tropical stable maps associated to $\barM_{0,n}(X, \beta)$ \cite{Ran}. Let $TSM_{\Gamma}(\Sigma_X, \beta) \subseteq TSM_{0,n}(\Sigma_X, \beta)$ be the subset for which the domain curve has dual graph $\Gamma.$ We have $TSM_{0,n}(\Sigma_X, \beta) = \bigcup_{\Gamma} TSM_{\Gamma}(\Sigma_X, \beta)$. 

For a dual graph $\Gamma$, consider $\barM_{\Gamma}(X, \beta) \subseteq \barM_{0,n}(X, \beta)$. We assume there exists a collection $\{U_{i, \Gamma}\}_i$ with classifying maps $b_{i, \Gamma}: U_{i, \Gamma} \rightarrow \barM_{\Gamma}(X, \beta)$ such that $\barM_{\Gamma}(X, \beta) = \bigcup_i b_{i, \Gamma}(U_{i, \Gamma})$, and $U_{i, \Gamma}$ gives rise to a product family of stable maps $\phi_{U_{i, \Gamma}}: (\mathcal{C}_{U_{i, \Gamma}} \times U_{i, \Gamma})/U_{i, \Gamma} \rightarrow X$ satisfying Assumption \ref{asu:locally_trivial}. Hence, we assume each $U_{i, \Gamma}$ is a toroidal, spin variety with spinor bundle $S_{U_{i, \Gamma}}$. We have $TSM_{\Gamma}(\Sigma_X, \beta) = \bigcup_{i} TSM_{\Gamma}(\Sigma_{X}, \beta)|_{b_{i, \Gamma}(U_{i, \Gamma})}$. For each $U_{i, \Gamma}$, we have maps $e^{K, trop}\left(\barN_{\left(\mathcal{C}_{U_{i, \Gamma}} \times U_{i, \Gamma}, \phi_{U_{i, \Gamma}}\right)}\right)^{-1}: A^{trop}_*(\Sigma_{U_{i, \Gamma}}) \rightarrow A^{trop}_*(\Sigma_{U_{i, \Gamma}}) \otimes \mathbb{C}(\kappa)$ from Definition \ref{def:tropical_Euler_local}. 

\begin{lemma}
\label{lem:strata}
For a dual graph strata $\barM_{\Gamma}(X, \beta) \subseteq \barM_{0,n}(X, \beta)$, there exists a well-defined map,

\[
e^{K, trop}\left(\barN_{n, \beta}|_{\barM_{\Gamma}(X, \beta)}\right)^{-1}: A^{trop}_*(TSM_{\Gamma}(\Sigma_X, \beta)) \rightarrow A^{trop}_* (TSM_{\Gamma}(\Sigma_X, \beta)) \otimes \mathbb{C}(\kappa)
\]
that computes the tropical, inverse Euler class of the SUSY bundle $\barN_{n, \beta}$ on $\barM_{\Gamma}(X, \beta)$.
\begin{proof}
    

    Suppose we have a collection $\{U_{i, \Gamma}\}_i$ of toroidal, spin varieties with classifying maps $b_{i, \Gamma}: U_{i, \Gamma} \rightarrow \barM_{\Gamma}(X, \beta)$ such that $\barM_{\Gamma}(X, \beta) = \bigcup_i b_{i, \Gamma}(U_{i, \Gamma})$. Let $\alpha \in A_*^{trop}(TSM_{\Gamma}(\Sigma_X, \beta))$. Consider the restriction \\$\alpha|_{Trop(b_{i, \Gamma})(U_{i, \Gamma})} \in A^{trop}_*(TSM_{\Gamma}(\Sigma_X, \beta)|_{Trop(b_{i, \Gamma})(\Sigma_{U_{i, \Gamma})}})$. Apply $e^{K, trop}\left(\barN_{\left(\mathcal{C}_{U_{i, \Gamma}} \times U_{i, \Gamma}, \phi_{U_{i, \Gamma}}\right)}\right)^{-1}$ to $(Trop(b_{i, \Gamma}))^*(\alpha|_{Trop(b_{i, \Gamma})(U_{i, \Gamma})}).$ This gives the tropical cycle,
    
    \[
    e^{K, trop}\left(\barN_{\left(\mathcal{C}_{U_{i, \Gamma}} \times U_{i, \Gamma}, \phi_{U_{i, \Gamma}}\right)}\right)^{-1}((Trop(b_{i, \Gamma}))^*(\alpha|_{Trop(b_{i, \Gamma})(U_{i, \Gamma})})) \in A_{*}^{trop}(\Sigma_{U_{i, \Gamma}})
    \]
    Pushforward with $(Trop(b_{i, \Gamma}))_*$ to get the cycle,
    
    \begin{align*}
       \alpha_i &:= (b^{trop}_{i,\Gamma})_*\left(e^{K, trop}\left(\barN_{\left(\mathcal{C}_{U_{i, \Gamma}} \times U_{i, \Gamma}, \phi_{U_{i, \Gamma}}\right)}\right)^{-1}\left(\left(Trop(b_{i, \Gamma})\right)^*\left(\alpha|_{Trop(b_{i, \Gamma})(U_{i, \Gamma})}\right)\right)\right) \\
       &\in A^{trop}_{*}\left(TSM_{\Gamma}(\Sigma_X, \beta)|_{b_{i,\Gamma}\left(U_{i, \Gamma}\right)}\right) \otimes \mathbb{C}(\kappa) 
    \end{align*}
    which is dependent on $U_{i, \Gamma}$. For $i \neq j$, the cycles $\alpha_i = \alpha_j$ will agree on intersections $b_{i, \Gamma}(U_{i, \Gamma}) = b_{j, \Gamma}(U_{j, \Gamma})$, because the restriction of $\barN_{n, \beta}$ onto $\barM_{\Gamma}$ is a vector bundle by Proposition \ref{prop:vector_bundle}, and by the functoriality of Chern classes with respect to pullback. Hence, we define $\alpha' := \cup_i \alpha_i \in A_*^{trop}(TSM_{\Gamma}(X, \beta)$. The map sending $\alpha \mapsto \alpha'$ is the desired map $e^{K, trop}\left(\barN_{n, \beta}|_{\barM_{\Gamma}(X, \beta)}\right)^{-1}$.
\end{proof}
\end{lemma}
Now, let $\alpha \in A^{trop}_*(TSM_{0,n}(\Sigma_X, \beta))$. For an $\barM_{\Gamma}(X, \beta) \subseteq \barM_{0,n}(X, \beta)$, consider the restriction $\alpha|_{TSM_{\Gamma}(\Sigma_X, \beta)} \in A^{trop}_*(TSM_{\Gamma}(\Sigma_X, \beta))$. We use the tropical, inverse Euler class map from Lemma \ref{lem:strata}. Since $\barN_{n, \beta}$ is a vector bundle on $\barM_{0,n}(X, \beta)$ by Proposition \ref{prop:vector_bundle} and that Chern classes are functorial under pullback, we have,

\[
e^{K, trop}\left(\barN_{n, \beta}|_{\barM_{\Gamma_i}(\Sigma_X, \beta)}\right)^{-1}\left(\alpha|_{TSM_{\Gamma_i}(\Sigma_X, \beta)}\right) = e^{K, trop}\left(\barN_{n, \beta}|_{\barM_{\Gamma_j}(X, \beta)}\right)^{-1}\left(\alpha|_{TSM_{\Gamma_j}(\Sigma_X, \beta)}\right)
\]
on $TSM_{\Gamma_i} \cap TSM_{\Gamma_j}$ for all $i,j$. Since $TSM_{0,n}(\Sigma_X, \beta) = \bigcup_{\Gamma} TSM_{\Gamma}$, we get a well-defined map that computes the tropical, inverse Euler class of $\barN_{n, \beta}$,

\begin{align*}
 e^{K, trop}(\barN_{n, \beta})^{-1}: A^{trop}_*(TSM_{0,n}(\Sigma_X, \beta)) &\rightarrow A^{trop}_*(TSM_{0,n}(\Sigma_X, \beta)) \otimes \mathbb{C}(\kappa) \\
 &\alpha \mapsto E^{trop}\left(\barN_{n, \beta}\right) \cap \alpha
\end{align*}
where $E^{trop}\left(\barN_{n, \beta}\right)$ is defined in Equation \ref{eq:inverse_Euler_exp_trop}. Thus, by the arguments made in Section \ref{sec:def_tropical_Euler},

\begin{theorem}
\label{thm:tropical_Euler_class}
There exists a well-defined map that computes the tropical, inverse Euler class of the SUSY normal bundle $\barN_{n, \beta} \rightarrow \barM_{0,n}(X, \beta)$,

\begin{align*}
   e^{K, trop}(\barN_{n, \beta})^{-1}: A^{trop}_*(TSM_{0,n}(\Sigma_X, \beta)) &\rightarrow A^{trop}_*(TSM_{0,n}(\Sigma_X, \beta)) \otimes \mathbb{C}(\kappa)\\
   \alpha &\mapsto E^{trop}\left(\barN_{n, \beta}\right) \cap \alpha 
\end{align*}
where $E^{trop}\left(\barN_{n, \beta}\right)$ is defined in Definition \ref{eq:inverse_Euler_exp_trop}.
\end{theorem}

\begin{remark}
\label{rem:upshot of the assumption}
    Instead of making Assumption \ref{asu:locally_trivial}, one could initially consider the splitting $\barN_{n, \beta} = \bigoplus_{i=1}^{rk_{n, \beta}}L_i$ as line bundles on $\barM_{0,n}(X, \beta)$, and consider $Trop(L_i)$ as tropical line bundles on $TSM_{0,n}(\Sigma_X, \beta)$. This is not ideal for a few reasons: 1) the definition of $\barN_{n, \beta}$ in Section \ref{sec:SUSY_bundles}, as well as the proof of Proposition \ref{prop:vector_bundle} requires the use of local neighborhoods $U \subseteq \barM_{0,n}(X, \beta)$, 2) tropicalization of line bundles \cite{CG} requires a toroidal base space 3) if we assume to have a cover $\{U_i\}_i$ of toroidal, spin varieties as in Assumption \ref{asu:locally_trivial}, this allows us to compute tropical Chern classes more concretely on the $U_i$, instead of working with all of $\barM_{0,n}(X, \beta)$. 
\end{remark}

\subsection{Tropical, super Gromov-Witten invariants}

Let $\left[\barM_{0,n}(X, \beta)\right]^{vir}$ be the virtual fundamental class \cite{BF}. Since $X$ is convex, $\left[\barM_{0,n}(X, \beta)\right]^{vir} = \left[\barM_{0,n}(X, \beta)\right]$. 

\begin{definition}
\label{def:tropical_SGW}
    The tropical, genus-0, $n$ Neveu-Schwarz marked, super Gromov-Witten invariant of a convex, toric variety $X$ in curve class $\beta \in H_2(X, \mathbb{Z})$ with incidence conditions $\gamma_i \in H^*(X)$ is,
    
    \begin{align*}
      SGW^{trop}_{0,n}(X, \beta) &:= e^{K, trop}(\barN_{n, \beta})^{-1}\left(Trop\left(\left[\barM_{0,n}(X, \beta)\right] \cap \prod_i ev_i^* \gamma_i\right)\right) \\
      &\in A^{trop}_*(TSM_{0,n}(\Sigma_X, \beta)) \otimes \mathbb{C}(\kappa)  
    \end{align*}
    where $e^{K, trop}(\barN_{n,\beta})^{-1}$ is defined in Theorem \ref{thm:tropical_Euler_class}.
\end{definition}

\begin{example}
\label{ex:SGW_P1}
Consider $\barM_{0,2}(\bP^1, H)$, the moduli space of genus-0, 2-marked stable maps to $\bP^1$ in the hyperplane class $H \in H_2(\bP^1, \mathbb{Z})$.  The moduli space is isomorphic to $\bP^1 \times \bP^1$, since the stable maps are determined by the images of the two marked points. It has virtual dimension 2, and the SUSY normal bundle $\barN_{2, H} \rightarrow \barM_{0,2}(\bP^1, H)$ is rank $rk_{2, H} = 2$. We show that,

\[
SGW^{trop}_{0,2}(\bP^1, H) :=  e^{K, trop}(\barN_{2, H})^{-1}\left(Trop\left(\left[\barM_{0,2}(\bP^1, H)\right]\right)\right)  = \frac{1}{\kappa^3}
\]
We recall from Section \ref{sec:tropicalization_lb} that $Trop\left(\barM_{0,2}(\bP^1, H)\right)$ is the constant function 1 on $TSM_{0,2}(\Sigma_{\bP^1}, H)$. To apply our tropicalization procedure, we verify that $\barM_{0,2}(\bP^1, H)$ satisfies Assumption \ref{asu:locally_trivial}. The moduli space has two dual graph strata: the smooth strata is isomorphic to $(\bP^1 \times \bP^1) \setminus \Delta \cong \cO_{\bP^1}(-2)$, with classifying map $b_1: \cO_{\bP^1}(-2) \rightarrow \barM_{0,2}(\bP^1, H)$. The other strata of stable maps with singular domain curve is isomorphic to the diagonal $\Delta \cong \bP^1$, with classifying map $b_2:\bP^1 \rightarrow \barM_{0,2}(\bP^1, H)$. Both strata are toric, and hence toroidal. Since $K_{\cO_{\bP^1}(-2)} \cong \cO_{\cO_{\bP^1}(-2)}$, both strata are also spin. The tropicalizations of the classifying maps, $Trop(b_1)$ and $Trop(b_2)$, are proper and give well-defined pullbacks on tropical Chow groups, since $b_1$ and $b_2$ are inclusions. 

We compute the tropical, inverse equivariant Euler class of $\barN_{2, H}$ on each strata. On $\cO_{\bP^1}(-2)$, the moduli space is given by the product family of stable maps $\phi_{\cO_{\bP^1}(-2)}: (\mathcal{C}_{\cO_{\bP^1}(-2)} \times \cO_{\bP^1}(-2))/\cO_{\bP^1}(-2) \rightarrow \bP^1$, where $\mathcal{C}_{\cO_{\bP^1}(-2)}$ is a family of smooth, genus-0, 2-marked curves. Therefore, we have $\barN_{2, H}|_{\cO_{\bP^1}(-2)} \cong \barN_{(\mathcal{C}_{\cO_{\bP^1}(-2)} \times \cO_{\bP^1}(-2), \phi_{\cO_{\bP^1}(-2)})} \cong \cO_{\cO_{\bP^1}(-2)}^{\oplus 2}$ by Example \ref{ex:product_SUSY}. Hence, $e^{K, trop}(\barN_{(\mathcal{C}_{\cO_{\bP^1}(-2)} \times \cO_{\bP^1}(-2), \phi_{\cO_{\bP^1}(-2)})})^{-1} = 0$.

Next, we compute the tropical, inverse equivariant Euler class of $\barN_{2, H}$ on $\Delta \cong \bP^1$. On $\Delta$, the moduli space of stable maps is given by the family $\phi_{\bP^1}: (\mathcal{C}_{\bP^1} \times \bP^1)/\bP^1 \rightarrow \bP^1$, where $C$ has two components, with 2 marked points on a contracted component, and the other component mapping to $H \in H_2(\bP^1, \mathbb{Z})$. The base $\Delta$ of $\phi_{\bP^1}$ parametrizes a marked point on $\bP^1$ that comes from the contracted component. Hence $\barN_{2, H}|_{\bP^1} \cong \barN_{(\mathcal{C}_{\bP^1}\times \bP^1, \phi_{\bP^1})}$. By Example \ref{ex:product_SUSY}, $N_{(\mathcal{C}_{\bP^1} \times \bP^1, \phi_{\bP^1})} \cong \cO_{\bP^1}(-1)^{\oplus 2}$. The tropical Chern class of $Trop(\cO_{\bP^1}(-1))$ is given in Example \ref{ex:chern_class_O(m)}. Hence, by Equation \ref{eq:inverse_Euler_exp_trop}, the contribution on $\Delta \cong \bP^1$ to $e^{K, trop}(N_{(\mathcal{C}_{\bP^1} \times \bP^1, \phi_{\bP^1})})^{-1}$ is $\frac{2}{\kappa^3}$. 

By \cite{Ful}, Proposition 1.8, the following sequence of Chow groups is exact for all $k$,

\[
A_k(\Delta) \xrightarrow{i_*} A_k(\bP^1 \times \bP^1) \xrightarrow{j^*} A_k(\cO_{\bP^1}(-2)) \rightarrow 0
\]
Consider $e(\barN_{2, H})^{-1} \in A_*(\bP^1 \times \bP^1)$. We showed that $j^*(e(\barN_{2, H})^{-1}) = e(\barN_{2, H}|_{\cO_{\bP^1}(-2)})^{-1} = e(\cO_{\cO_{\bP^1}(-2)}^{\oplus 2})^{-1} = 0 \in A_*(\cO_{\bP^1}(-2))$. By exactness, there exists $\alpha \in A_*(\Delta)$ such that $i_* \alpha = e(\barN_{2, H})^{-1}$. By the self-intersection formula, $i^*i_* \alpha = e(N) \cdot \alpha = e(\barN_{2, H}|_{\Delta})^{-1}$, where $N$ is the normal bundle of $\Delta$, hence $N \cong T\bP^1 \cong O_{\bP^1}(2)$. Since $e(N)$ is invertible in the Chow ring, this implies $\alpha = e(\barN_{2, H}|_{\Delta})^{-1} \cdot e(N)^{-1}$. We have $\deg(e(\barN_{2, H})^{-1}) = \deg(i_*\alpha) = \deg(\alpha)$ by \cite{Ful}, Definition 1.4. By the above computation, we have $\deg \alpha = -2\cdot\frac{1}{2} = -1$. By Equation \ref{eq:inverse_Euler_exp}, this implies that $\deg (e^K(\barN_{2, H}|_{\Delta})^{-1} \cdot e(N)^{-1}) = \frac{1}{\kappa^3}$. Since Chern classes are analogous to tropical Chern classes for $\bP^1$ by Example \ref{ex:chern_class_O(m)}, we have $\deg (e^{K, trop}(\barN_{2, H}|_{\Delta})^{-1} \cdot e(N)^{-1}) = \frac{1}{\kappa^3}$. Therefore, $\deg(e^{K, trop}(\barN_{2, H})^{-1}) = \frac{1}{\kappa^3}$. Thus, since there are also no incidence conditions $\gamma_i$, we have $SGW_{0,2}^{trop}(\bP^1, H) = \frac{1}{\kappa^3}$. 
\end{example}

\begin{corollary}
\label{cor:connect}
    Let $X = pt$. By evaluating the tropical, genus-0, $n$-marked, super Gromov-Witten invariant defined in Definition \ref{def:tropical_SGW} over the smooth locus $\left[\mathcal{M}_{0,n}\right] \subset \left[\barM_{0,n}\right]$, we recover Theorem \ref{thm:TSGW_point}.

    \begin{proof}
        Recall that when $X = pt$, we have $\barM_{0,n}(pt, 0) \cong \barM_{0,n}$, and we have the isomorphism of SUSY bundles $\barN_{n, \beta} \cong \barN_n$. By \cite{CG}, Corollary 6.22, the description of tropical Chern classes from \cite{CGM} recovers the description of tropical $\psi$-classes in Definition \ref{def:tropical_psi}. Since Theorem \ref{thm:tropical_Euler_class} relies on \cite{CGM}, the tropical, inverse Euler class of $\barN_{n, \beta} \cong \barN_n$ will recover the tropical descendant invariants considered in \cite{KM}. Thus, for $X = pt$, we recover Theorem \ref{thm:TSGW_point} by evaluating the tropical, inverse Euler class over $\left[\mathcal{M}_{0,n}\right] \subset \left[\barM_{0,n}\right]$.
    \end{proof}
\end{corollary}

\bibliographystyle{alpha}
\bibliography{refs}

\end{document}